\input amssym.def 
\input amssym
\magnification=1200
\parindent0pt
\hsize=16 true cm
\baselineskip=13  pt plus .2pt
$ $

\def\Z{{\Bbb Z}}
\def\D{{\Bbb D}}
\def\A{{\Bbb A}}
\def\S{{\Bbb S}}

\centerline {\bf  A note on large bounding and non-bounding finite group-actions}

\centerline {\bf on surfaces of small genus}

\bigskip

\centerline {Bruno P. Zimmermann}

\medskip

\centerline {Universit\`a degli Studi di Trieste}

\centerline {Dipartimento di Matematica e Geoscienze}

\centerline {34127 Trieste, Italy}

\bigskip \bigskip

{\bf Abstract.}  The classification of finite group-actions on closed 
surfaces of small genus is well-known. In the present paper we are interested
in the question of which of these group-actions are bounding (extend to a
compact 3-manifold with the surface as its unique boundary component, e.g. to a
handlebody) or geometrically bounding  (extend to a hyperbolic 3-manifold with
totally geodesic boundary)
concentrating, as a typical case, on large group-actions on surfaces of genus 3.

\bigskip \bigskip

{\bf 1. Introduction}

\medskip

All finite group-actions in the present paper will be
orientation-preserving and faithful, all manifolds will be orientable.
We are interested in large group-actions  of a finite group $G$ on a closed
hyperbolic surface
$\Sigma$ of genus $g \ge 2$.  By choosing a hyperbolic structure on a 
quotient-orbifold
$\Sigma/G$ and lifting the structure to $\Sigma$, we can assume that $G$ acts by
isometries, for some hyperbolic structure on $\Sigma$; then  the group of all lifts
of elements of
$G$ to the universal covering $\Bbb H^2$ of $\Sigma$ is a Fuchsian group $F$, and we
have an exact sequence
$$1 \to K \hookrightarrow  F  \to G \to 1$$
where $K \cong \pi_1(\Sigma)$ denotes the universal covering group. 
We denote a  Fuchsian group by its signature (which is also the signature of the
quotient-orbifold $\Sigma/G$); for example, by $(p,q,r)$ we denote the triangle
group (of genus 0 which we omit) of orientation-preserving elements in the group
generated by the reflections in the sides of a hyperbolic triangle with angles
$2\pi/p$, $2\pi/q$ and $2\pi/r$, and a quadrangle group $(p,q,r,s)$ is defined
analogously.

\medskip

We say that a finite $G$-action on a surface $\Sigma$ {\it bounds} if the
$G$-action on $\Sigma$ extends to a $G$-action on a compact 3-manifold $M$ with 
$\partial M  = \Sigma$ (e.g., to a handlebody); it {\it bounds geometrically} if $M$
can be chosen as a  compact hyperbolic 3-manifold $M$ with totally geodesic
boundary. In this second case, by an application of Mostow rigidity one can assume
that $G$ acts by isometries also on $M$ (and then also $\Sigma = \partial M$
achieves a hyperbolic structure on which $G$ acts by isometries).

\medskip

In [WZ] the authors determine which
finite group-actions on a surface of genus two extend to a 3-manifold (i.e.,
bound),  and in particular also to the 3-sphere, for some embedding of the surface
into $S^3$ (in order to give explicit geometric descriptions or "visualizations
of these actions in the familiar 3-space). In the present paper we consider the case
$g = 3$ but concentrate  on possible extensions to handlebodies and hyperbolic
3-manifolds with totally geodesic boundary instead of $S^3$.

\medskip

We refer to Broughton [B] for the classification of the finite group-actions on
surfaces of genus 3.  We will represent
finite group-actions on surfaces by surjections $F \to G$, always assumed to have
torsionfre kernel, of a Fuchsian group $F$ onto a finite group $G$; the kernel of
such a surjection is a torsionfree Fuchsian group $K$ (a surface group), and $G
\cong F/K$ acts (by isometries) on the hyperbolic surface $\Bbb H^2/K$.
In the following theorem, we list the largest group-actions on a surface of genus 3
in decreasing order and determine which of these actions bound, bound a handlebody
or bound geometrically, representing the actions by a surjection $F \to G$ of a
Fuchsian group $F$ to a finite group $G$. We denote by $\D_n$ the dihedral group of
order $2n$, by $\A_4$ and $S_4$ the alternating and symmetric groups of orders 12
and 24; for the groups $D_{2,8,5}$ and $D_{2,12,5}$ we refer to [B].

\bigskip

{\bf Theorem.}   {\sl  The bounding and non-bounding finite group-actions on
a surface of genus 3, of order $\ge 24$, are the following.

\medskip
 
i) The two largest group-actions, 
of orders 168 and 96 are  represented by surjections
$$(2,3,7) \to  {\rm PSL}_2(7) \;\; {\sl and} \;\; (2,3,8)  \to  
\D_3 \ltimes (\Z_4 \times \Z_4)$$

and do not bound. 

\medskip

ii) Two actions of order 48 associated to surjections
$$(3,3,4)  \to \Z_3 \ltimes (\Z_4 \times \Z_4) \;\; {\it and} \;\;
(2,4,6) \to  \Z_2 \times S_4;$$ 

the first one is a subgroup of index 2 of the group of order 96 in i) and does not
bound, the second one is the largest bounding group-action on a surface of genus 3;
it bounds geometrically but does not extend to  a handlebody.

\medskip

iii) Two non-bounding actions of order 32 associated to surjections
$$(2,4,8)  \to \Z_2 \ltimes (\Z_2 \times \Z_8) \;\; {\it and} \;\;
(2,4,8)  \to \Z_2 \times D_{2,8,5}.$$

\medskip

iv) Two non-bounding actions of order 24  associated to surjections

$$(3,3,6) \to  {\rm SL}_2(3)  \;\; {\it and} \;\;  (2,4,12) \to  D_{2,12,5}.$$

\medskip
\vfill  \eject

v)  Three bounding actions of order 24 associated to surjections

$$(2,6,6) \to \Z_2 \times \A_4, \;\;  (3,4,4) \to \S_4  \;\; {\it and} \;\; 
(2,2,2,3) \to \S_4;$$ 

these three actions are 
subgroups of index 2 of the geometrically bounding action of order 48 in ii). The
last one is also the largest action on a surface of genus 3 which extends to a
handlebody (in fact $\S_4$ is the unique maximal handlebody group of genus 3,  
of maximal possible order $12(g-1)$).

\medskip

Finally, for each of the non-bounding actions in i) - iv) there is already a
cyclic subgroup which does not bound. }

\bigskip

{\bf Corollary.}   {\sl The largest bounding finite group-action on a surface of
genus 3 is an action of $\Z_2 \times S_4$ of order 48 which bounds geometrically
but does not extend to a handlebody. The largest finite group-action in genus 3
which extends to a handlebody is the action of the subgroup $S_4$ which bounds also
geometrically.}

\bigskip

The largest finite group-action on a surface of genus 4 is an action of the
symmetric group $\S_5$ associated to a surjection $(2,4,5)  \to \S_5$ 
(cf. [C] and its references); using similar methods, the following holds.

\bigskip

{\bf Proposition.}   {\sl The largest group-action on a surface of genus 4, of type 
$(2,4,5)  \to \S_5$, bounds geometrically. The largest action in genus 4 which
extends to a handlebody is of type   $(2,2,2,3)  \to \D_3 \times \D_3$.}

\bigskip

See also [Z4] for a discussion of various aspects of finite group-actions on
surfaces.

\bigskip

{\bf 2. Proof of the Theorem}

\medskip

i)  The Hurwitz-action of
${\rm PSL}_2(7)$ on Klein's quartic $\Sigma_3$ of genus 3 does not bound, i.e. does 
not extend to any compact 3-manifold
$M$ with exactly one boundary component $\partial M = \Sigma_3$: the quotient
orbifold  $\Sigma_3/{\rm PSL}_2(7)$ is the 2-sphere with three branch points of
orders 2, 3 and 7 which does not occur as the unique boundary component of a
compact 3-orbifold since a singular axis starting in the boundary point of order 7
can end only in a dihedral point of dihedral type $\D_7$ but ${\rm PSL}_2(7)$ has no
dihedral subgroup $\D_7$. A cyclic subgroup which does not bound is of type 
$(7,7,7) \to \Z_7$.

\medskip

Similarly, the group $\D_3 \ltimes (\Z_4 \times \Z_4)$ acting on Fermat's quadric of
genus 3 has no subgroup $\D_8$ and hence does not bound, a cyclic non-bounding
subgroup is of type  $(4,8,8) \to  \Z_8$.

\medskip

ii)  Concerning the first case, an axis of order 4 starting in a singular
point of order 4 in the quotient 2-orbifold of type (3,3,4) can only end in a
singular point of type $\D_4$ or $\S_4$ but  $\Z_3 \ltimes (\Z_4 \times \Z_4)$ has
no such subgroups;  a cyclic non-bounding
subgroup is of type $(4,4,4,4) \to  \Z_4$.

\medskip

We show that the action of $\Z_2 \times S_4$ bounds geometrically. We consider 
a hyperbolic tetrahedron $\cal T$, truncated by an orthogonal hyperplane at a 
vertex of type  (2,4,6) where three edges of orders 2, 4 and 6 meet (with
dihedral angles $\pi/2, \pi/4$ and $\pi/6$); the edge opposite to the edge of
singular order 6 has order 3, all other edges have singular order 2. Let $T$ denote
the associated tetrahedral group (the orientation-preserving subgroup of index
2 in the Coxeter group generated by the reflections in the four faces of the
tetrahedron), with a triangle subgroup of type (2,4,6) generated by the rotations 
at the truncated vertex. Then it is easy to see that a surjection  
$(2,4,6) \to  \Z_2 \times S_4$
defining the action on the surface of genus 3 extends to a surjection $T \to  \Z_2
\times S_4$. The covering of the 3-orbifold $\cal T$ associated to the kernel of
this surjection is a hyperbolic 3-manifold with totally geodesic boundary (a
surface of genus 3), with an isometric action of $\Z_2 \times S_4$; this restricts
to an isometric action of  $\Z_2 \times S_4$ on the boundary which realizes the
action in part ii) of the Theorem.

More explicitly, the rotational generators of the tetrahedral group $T$,
suitably oriented, can be mapped to the following elements of $\Z_2 \times S_4$ 
where $c$ denotes the generator of $\Bbb Z_2$:
the order 4 rotation is mapped to  $c \, (1234)$, the order 6 rotation to $c \,
(143)$ and the order 3 rotation to (142), and this determines also the images of 
the order 2 rotations (see [GZ] or [Z2] for similar constructions).

\medskip

The action of $\Z_2 \times S_4$ does not extend to a handlebody. More generally,
any action $(p,q,r)  \to G$  associated to a triangle group $(p,q,r)$ does not
extend to a handlebody, see [Z1].

\medskip

iii) and iv) are similar to i). The cyclic non-bounding subgroups in iii) are both
of type $(4,8,8) \to  \Z_8$, in iv) of types $(2,3,3,6) \to  \Z_6$  and
$(2,12,12) \to  \Z_{12}$.  Concerning v), it is well-known that $\S_4$ is the
unique maximal handlebody group of genus 3 (of order $12(g-1) = 24$), see [Z3].

\bigskip

The {\it Proof of the Proposition} is similar to the geometrically
bounding case  of part ii) of the Theorem, with the following choices.  We consider 
a hyperbolic tetrahedron, truncated by an orthogonal hyperplane at a 
vertex of type  (2,4,5) where three edges of orders 2, 4 and 5 meet; an 
edge of order 3 connects the edges of orders 4 and 5, all other edges have
order 2.  The order 4 rotation is mapped to  $(2345)$, the order 5 rotation to
$(12345)$ (oriented such that (12345) (12) = (2345)), the order 3 rotation to
(135) (such that (12)(34) (135) = (12345)). With these choices, the proof is the
same as that of part ii) of the Theorem.

\medskip

Finally, $\D_3 \times \D_3$ is the unique
maximal handlebody group of genus 4, of maximal possible order $12(g-1)$ (see
[Z3]).

\bigskip  \bigskip
\vfill   \eject

\centerline {\bf References}

\medskip

\item {[B]} S.A. Broughton, {\it Classifuing finite group actions on
surfaces of low genus.}  J. Pure Appl. Alg. 69 (1990),  233-270

\item {[C]} M.D.E. Conder, {\it Large group actions on surfaces.} Contemp. Math.
629  (2014),  77-98

\item {[GZ]} M. Gradolato, B. Zimmermann,  {\it  Extending finite group
actions on  surfaces to hyperbolic 3-manifolds.} Math. Proc. Cambridge
Phil. Soc.  117   (1995), 137-151

\item {[WZ]}  C. Wang, S. Wang, Y. Zhang, B. Zimmermann,  {\it  Finite group
actions on the genus-2 surface, geometric generators and extendability.} 
Rend. Istit. Mat. Univ. Trieste 52  (2020), 513-524  (electronic version under
http://rendiconti.dmi.units.it)

\item {[Z1]} B. Zimmermann,  {\it \"Uber Abbildungsklassen von
Henkelk\"orpern.}  Arch. Math. 33 (1979), 379-382

\item {[Z2]}  B. Zimmermann, {\it  Hurwitz groups and finite group actions on
hyperbolic 3-manifolds.}  J. London Math. Soc. 52  (1995), 199-208

\item {[Z3]} B. Zimmermann,  {\it Genus actions of finite groups on
3-manifolds.}  Michigan Math. J. 43 (1996),  593-610

\item {[Z4]} B. Zimmermann,  {\it Hurwitz groups, maximal reducible groups and
maximal handlebody groups.}  arXiv:2110.11050

\bye